\input amstex
\magnification=\magstep1 
\baselineskip=13pt
\documentstyle{amsppt}
\CenteredTagsOnSplits \NoRunningHeads
\vsize=8.7truein 
\def\today{\ifcase\month\or
  January\or February\or March\or April\or May\or June\or
  July\or August\or September\or October\or November\or December\fi
  \space\number\day, \number\year}

\def\PP{{\bold P}}
\def\EE{{\bold E\thinspace }}

\def\per{\operatorname{per}}
\def\Mat{\operatorname{Mat}}

\topmatter
\title Counting magic squares in quasi-polynomial time \endtitle
\author Alexander Barvinok, Alex Samorodnitsky, and  Alexander Yong\endauthor
\address Department of Mathematics, University of Michigan, Ann Arbor,
MI 48109-1043, USA \endaddress
\email barvinok$\@$umich.edu \endemail
\thanks The research of the first author was partially supported by NSF Grant DMS 0400617.
The research of the second author was partially supported by ISF grant 039-7165.
The research of the third author was partially completed while he  was an NSF sponsored visitor
at the Institute for Pure and Applied Mathematics at UCLA, during April-June 2006. 
The third author was also partially supported by NSF grant 0601010 and an NSERC Postdoctoral
fellowship held at the Fields Institute, Toronto. 
\endthanks
\address Department of Computer Science, Hebrew University of Jerusalem, Givat Ram Campus,
91904, Israel \endaddress
\email salex$\@$cs.huji.ac.il \endemail
\address Department of Mathematics, University of Minnesota, Minneapolis, MN 55455, USA
\endaddress
\email ayong$\@$math.umn.edu \endemail
\date March 2007 \enddate
\abstract  We present a randomized algorithm, which, given positive integers $n$ and $t$
and a real number $0< \epsilon <1$, computes the number $|\Sigma(n, t)|$ of $n \times n$ non-negative
integer matrices (magic squares)  with the row and column sums equal to $t$ within relative error 
$\epsilon$. The computational complexity of the algorithm is polynomial in $\epsilon^{-1}$ and quasi-polynomial in $N=nt$, that is, of the order $N^{\log N}$.  A simplified version of the algorithm works in 
time polynomial in $\epsilon^{-1}$ and $N$ and estimates $|\Sigma(n,t)|$ within a factor of 
$N^{\log N}$. This simplified version has been implemented. We present results of the implementation,
state some conjectures, and discuss possible generalizations. \endabstract
\endtopmatter
\document

\head 1. Introduction and main results \endhead

\subhead (1.1) Magic squares \endsubhead
Let us fix two positive integers $n$ and $t$. An $n \times n$ {\it magic square with the line sum} $t$
is an $n \times n$ non-negative integer matrix $D=\left(d_{ij}\right)$ with the row and column
sums $t$:
$$\split &\sum_{j=1}^n d_{ij}=t \quad \text{for} \quad i=1, \ldots, n \quad \text{and} \\
&\sum_{i=1}^n d_{ij}=t \quad \text{for}\quad j=1, \ldots, n. \endsplit$$
We note that sometimes such matrices are called semi-magic squares, but we follow the
terminology adopted in modern combinatorics, for example, in Chapter 4 of \cite{St97}. 

Let $\Sigma(n, t)$ be the set of all $n \times n$ magic squares with the line sum $t$. In this 
paper, we present a randomized approximation algorithm to compute the number 
$|\Sigma(n,t)|$. The algorithm runs in quasi-polynomial time. More precisely, let 
$N=nt$ (in what follows, we reserve notation $N$ for the sum of the entries of the matrix).
 \medskip
  We present a randomized algorithm, which, for any given positive integers $n$ and $t$ and positive $\epsilon <1$, approximates $|\Sigma(n,t)|$ within relative error $\epsilon$. The computational complexity of the algorithm is
$ \left(1/\epsilon \right)^{O(1)} N^{O(\ln N)}$
(in the unit cost model).
\medskip

From this same approach, one also obtains a simpler, randomized polynomial time algorithm
which approximates  $|\Sigma(n,t)|$ within a factor of $N^{O(\ln N)}$.
We implemented the latter algorithm and report on the computational results in Section 1.4.

\subhead (1.2) Contingency tables \endsubhead
More generally, given positive integers $m$ and $n$, a positive integer vector 
$R=\left(r_1, \ldots, r_m \right)$, and a positive integer vector $C=\left(c_1, \ldots, c_n \right)$
such that 
$$r_1 +\ldots +r_m=c_1 + \ldots + c_n =N,$$
the $m \times n$ non-negative integer matrices with the row sums $r_1, \ldots, r_m$ and the 
column sums $c_1, \ldots, c_n$ are called {\it contingency tables} with the margins $R$ and $C$. The problem of computing or estimating efficiently the cardinality $|\Sigma(R,C)|$ of the set 
of contingency tables with the given margins has been of significant interest, see \cite{DE85},
\cite{DG95},
\cite{D+97}, \cite{Mo02}, and \cite{CD03} due to  connections to statistics, representation theory,
and symmetric functions.

Using the Markov Chain Monte Carlo approach, Dyer, Kannan, and Mount \cite{D+97} showed how to count contingency tables if the row and column sums are sufficiently
large, namely, if $r_i=\Omega\left(n^2 m\right)$ and $c_j=\Omega\left(m^2 n \right)$ for all $i,j$. 
They presented a randomized algorithm, which, given an $\epsilon>0$, approximates the 
number $|\Sigma(R,C)|$ of tables within relative error $\epsilon$ in time polynomial in
$\epsilon^{-1}$, $n$, $m$, and $\sum_i \log r_i + \sum_j \log c_j$ (the bit size of the margins).
It turns out that for large margins the number of contingency tables is well-approximated by the
volume of the {\it transportation polytope} of the $m \times n$ non-negative matrices with 
the row sums $r_i$ and the column sums $c_j$. The set $\Sigma(R,C)$ of contingency 
tables can be viewed as the set of integer points in that polytope. Subsequently, Morris \cite{Mo02}
obtained a similar result for the bounds $r_i=\Omega\left(n^{2/3} m \ln m \right)$ and $c_j=\Omega\left(m^{3/2} n \ln n\right)$. 

In addition,  for large $r_i$ and $c_j$ there is a heuristic formula for $|\Sigma(R,C)|$ 
due to Diaconis and Efron \cite{DE85}:
$$\split |\Sigma(R, C)| \approx &\left({2N+ mn \over 2} \right)^{(m-1)(n-1)}  \left(\prod_{i=1}^m \overline{r_i} \right)^{n-1} \\ &\times \left( \prod_{j=1}^n \overline{c_j} \right)^{k-1}  {\Gamma(nk) \over \Gamma^m(n)\Gamma^n(k)},  \endsplit \tag1.2.1$$
where
$$\split &\overline{r_i}={1-w \over m} +{w r_i \over N} \quad \text{and} \quad c_j={1- w \over n} + {w c_j \over N}  \\ &\text{for} \quad w={1 \over 1 +mn/2N} 
\quad \text{and} \quad k={n+1 \over n \sum_{i=1}^m \overline{r_i}^2} -{1 \over n}. \endsplit$$
This formula first approximates $|\Sigma(R,C)|$ by the volume of the corresponding transportation 
polytope and then (since no explicit formula is known for the volume) approximates the volume 
by an explicitly computable integral of a certain density. However, there are no proved or even 
conjectural conclusions on the accuracy of this formula.

At the opposite extreme, when the margins $r_i, c_j$ are very small (bounded by a constant fixed in advance) relative to the dimensions $m$ and $n$ of the matrix,  B\' ek\' essy, B\' ek\' essy, and Koml\' os \cite{B+72} obtained an asymptotic formula
$$|\Sigma(R,C)| \approx {N! \over r_1! \cdots r_m! c_1! \cdots c_n!} \exp\left\{{2 \over N^2}
\sum_{i,j} {r_i \choose 2} {c_j \choose 2} \right\}. \tag1.2.2$$
This formula reflects the fact that the majority of contingency tables with small 
margins have entries $0$, $1$, and $2$.
Also, if the margins are bounded by a constant fixed in advance, one can compute the exact value of $|\Sigma(R,C)|$ in time 
polynomial in $m+n$ by a dynamic programming algorithm.

Using the dynamic programming approach in a different vein, Cryan and Dyer \cite{CD03} constructed a randomized polynomial time approximation
algorithm for computing $|\Sigma(R,C)|$ provided the number $m$ of rows (or the number $n$ 
of columns) is fixed in advance.

In some sense, the case of magic squares $m=n$ with moderately large margins $t$ (say, of 
order $n$) lies at the core of the remaining hard cases of contingency table enumeration. Our 
algorithm is quasi-polynomial, and we conjecture that its straightforward modification achieves,
in fact, a genuinely polynomial time complexity, and that it naturally extends to a randomized polynomial
time algorithm to count contingency tables with any margins, see Section 10. A disadvantage of our approach
is that there seems to be no easy way to generate a {\it random} magic square, unlike 
in the approaches of \cite{D+97} and \cite{Mo02}.

\subhead (1.3) Idea of the algorithm \endsubhead Our algorithm builds on the 
technique of rapidly mixing Markov chains and, in, particular on efficient sampling from log-concave densities, as developed in \cite{AK91}, \cite{F+94}, \cite{FK99}, \cite{LV06}, see also \cite{Ve05} for a survey, the permanent approximation algorithm
\cite{J+04}, the strongly polynomial time algorithm for matrix scaling \cite{L+00}, as well as the integral representation of $|\Sigma(R,C)|$ from \cite{Ba05} and 
\cite{Ba07}.

Let $\Delta=\Delta_{n \times n} \subset {\Bbb R}^{n^2}$ be the open $(n^2-1)$-dimensional simplex of all $n \times n$ positive matrices $X=\left(x_{ij}\right)$
such that 
$$\sum_{i,j=1}^n x_{ij}=1.$$
 Let $\mu$ be Lebesgue measure
on $\Delta$ normalized by the constraint $\mu(\Delta)=1$.
 Using results of \cite{Ba05}, we 
represent the number of magic squares as an integral
$$|\Sigma(n,t)|=\int_{\Delta} f \ d \mu \tag1.3.1$$
of some continuous density $f: \Delta \longrightarrow {\Bbb R}_+$. Furthermore, as in \cite{Ba05},
the density $f$ is factored 
$$f=p \phi, \tag1.3.2$$
where $\phi: \Delta \longrightarrow {\Bbb R}_+$ is {\it log-concave}, that is,
$$\split \phi(\alpha X +\beta Y) \geq \phi^{\alpha}(X) \phi^{\beta}(Y) \quad& \text{for all} \quad 
X, Y \in \Delta \quad \text{and} \\ & \text{for all} \quad \alpha,\beta \geq 0 \quad \text{such that} \quad 
\alpha+\beta=1 \endsplit$$
and $p(X) \geq 1$ for all $X \in \Delta$.
Moreover, for any $X \in \Delta$ the values $p(X)$ and $\phi(X)$ are computable in time polynomial
in $N$. More precisely, for any given $\epsilon>0$ the value of $\phi$ can be computed 
within relative error $\epsilon$ in time polynomial in $\ln (1/\epsilon)$ and $N$ by a deterministic 
algorithm of \cite{L+00} while the value of $p$ can be computed within relative error $\epsilon$ in time polynomial
in $1/\epsilon$ and $N$ by a randomized algorithm of \cite{J+04}. The algorithm computing $\phi$ 
seems to work very well in practice.

The key result of this paper is that there is a {\it threshold} $T=N^{\kappa \ln N}$ for some absolute 
constant $\kappa>0$ such that 
if we define $\overline{p}: \Delta \longrightarrow {\Bbb R}_+$ by 
$$\overline{p}(X)=\cases p(X) &\text{if \ } p(X) \leq T \\ T &\text{if \ } p(X) >T \endcases$$
then the integral
$$\int_{\Delta} \overline{p} \phi \ d\mu \tag1.3.3$$
approximates the integral
$$\int_{\Delta} p \phi \ d \mu =|\Sigma(n,t)|$$
within a relative error as small as $N^{-n}$, say. 

The simplified version of the algorithm consists of computing the integral 
$$\int_{\Delta} \phi \ d \mu \tag1.3.4$$
using any of the randomized polynomial time algorithms of  \cite{AK91}, \cite{F+94}, \cite{FK99}, and \cite{LV06}
for {\it integrating} log-concave densities. Clearly, the value of (1.3.4) approximates $|\Sigma(n,t)|$
within a factor of $N^{O(\ln N)}$, that is,
$$\int_{\Delta} \phi \ d \mu \leq |\Sigma(n, t)| \leq N^{\kappa \ln N} \int_{\Delta} \phi \ d \mu$$
for some absolute constant $\kappa \geq 0$.

The full version of the algorithm consists of estimating (1.3.3) within relative error 
$\epsilon$ using any of the randomized polynomial time algorithms  of \cite{AK91}, \cite{F+94}, \cite{FK99}, and \cite{LV06} for {\it sampling}
from log-concave densities as well as the randomized polynomial time algorithm of \cite{J+04}
for approximating the permanent of a positive matrix. 

Namely, let $\nu$ be the Borel probability measure on $\Delta$ with the 
density proportional to $\phi$.

Thus we can rewrite (1.3.3) as the product
$$\left(\int_{\Delta} \overline{p} \ d \nu  \right) \left( \int_{\Delta} \phi \  d \mu \right).$$
We compute the second factor as above.
The first factor is approximated by the sample mean
$$\int_{\Delta} \overline{p} \ d \nu \approx {1 \over m} \sum_{i=1}^m \overline{p}(x_i), \tag1.3.5$$
where $x_1, \ldots, x_m \in \Delta$ are independent points sampled at random from the measure $\nu$.
Each such point can be sampled in time polynomial in $N$.
The Chebyshev inequality implies that to achieve relative error $\epsilon$
with probability $2/3$ it suffices to sample $m=O\left(\epsilon^{-2} T^2 \right)=\epsilon^{-2} N^{O(\ln N)}$
points in (1.3.5).

\subhead (1.4) Computational experiments \endsubhead 

We implemented the simplified version of the algorithm that computes the
integral (1.3.4). Below, we have tabulated some representative examples
of our results for various values of
$t$ and $n$. The software and additional data are available at \cite{Yo07}.

For $n \leq 7$, we were able to compare the obtained values (``estimate'' in the table below)
against the exact numbers (``answer'' in the table below is the exact number rounded to the 
three main digits) computed at our request by Jes\'us De Loera using the {\tt LattE} code,
see \cite{L+04}:

$$\matrix &n &\quad  &t  &\quad  &\text{estimate} &\quad &\text{answer}\\
&\quad &\quad &\quad &\quad &\quad &\quad &\quad \\
                &5  &\quad &5  &\quad &1.51\times 10^7 &\quad  &2.21\times 10^7  \\
                &\phantom{5}  &\quad &10 &\quad &7.96 \times 10^{10}   &\quad &7.93\times 10^{10} \\
                &\phantom{5} &\quad &125  &\quad   &1.55\times 10^{27}  &\quad  &1.10 \times 10^{27}\\                
&\quad &\quad &\quad &\quad &\quad &\quad &\quad \\
&6 &\quad   &6  &\quad    &4.70 \times 10^{11} &\quad     &6.02 \times 10^{11}\\
 &\phantom{6} &\quad    &12 &\quad    & 7.12\times 10^{16} &\quad    &2.28\times 10^{17} \\
  &\phantom{6} &\quad   & 36  &\quad     &2.21\times 10^{27}  &\quad    &5.62 \times 10^{27} \\
   &\phantom{6} &\quad    &216   &\quad    & 3.76 \times 10^{46} &\quad     &3.07 \times 10^{46} \\
&\quad &\quad &\quad &\quad &\quad &\quad &\quad \\ 
 &7  &\quad     &7    &\quad       &1.70\times 10^{17} &\quad   &   2.16\times 10^{17} \\
&\phantom{7}  &\quad  &14  &\quad     &1.75\times 10^{25}  &\quad    &2.46 \times 10^{25} \\
 &\phantom{7} &\quad   & 49  &\quad     & 3.00\times 10^{42} &\quad   & 4.00\times 10^{42} \\
  &\phantom{7} &\quad    &343  &\quad    & 3.90\times 10^{71} &\quad   &1.28\times 10^{72}.
 \endmatrix$$
 
The larger values of $n$ below seem to be presently beyond the reach of {\tt LattE} or any other available computer code. Thus we compared the obtained values (``estimate'' in the table below) with the Diaconis-Efron 
heuristic formula (1.2.1)  (``heuristic'' in the table below), which is believed to be valid for $t \gg n$:
                 
$$\matrix &n &\quad &t &\quad &\text{estimate} &\quad &\text{heuristic} \\
&\quad &\quad &\quad &\quad &\quad &\quad &\quad \\
&12 &\quad  &8 &\quad & 2.72 \times 10^{48} & \quad &4.96 \times 10^{49} \\
&\phantom{12} & \quad & 20 & \quad & 4.55 \times10^{81} & \quad & 1.68 \times 10^{82} \\
&\quad &\quad &\quad &\quad &\quad &\quad &\quad \\
&15 &\quad &20 &\quad & 8.3 \times 10^{119} & \quad &2.43 \times 10^{121} \\
&\phantom{15} &\quad &100 &\quad &7.65 \times 10^{236} &\quad &2.71 \times 10^{237} .\\ 
\endmatrix$$

An interesting feature of the data is that the Diaconis-Efron formula seems to be in a reasonable 
agreement with our computations even at the range $t \sim n$, where the heuristic arguments supporting the formula do not look plausible any longer. One can argue, however,  that the agreement
becomes more reasonable for larger $t$.

Finally, we consider the case of very small $t$, where we compare our results (``estimate'' in the 
table below) with the asymptotic formula (1.2.2)  (``asymptotic'' in the table below).

$$\matrix &n &\quad &t &\quad &\text{estimate} &\quad &\text{asymptotic} \\
&\quad &\quad &\quad &\quad &\quad &\quad &\quad \\
&25 &\quad &5 &\quad &2.89 \times 10^{108} &\quad &6.17 \times 10^{108} \\ 
&30 &\quad &5 &\quad & 1.49 \times 10^{142} & \quad &3.02 \times 10^{142}. \endmatrix$$

These results suggest that the integral (1.3.4) approximates $|\Sigma(n, t)|$
quite well, conjecturally within a factor of $N^{O(1)}$, if not just $O(1)$. Moreover,
the algorithm appears to estimate $|\Sigma(R,C)|$ for general
margins $(R,C)$ with similar effectiveness. For example, in a
known test case \cite{DE85}, \cite{L+04}, we have $n=4$,
$R=(220,215,93,64)$ and $C=(108,286,71,127)$.
The correct answer is about $1.22\times 10^{15}$ whereas
our algorithm predicts it to be $1.02\times 10^{15}$.

In practice, our algorithm is slower than {\tt LattE} when the latter is
applicable. However, as is explained in Section 3.6,
our algorithm has the advantage of being able to compute
when the values of $n$ are moderately large, far beyond those possible with exact
methods. Typically, we have been able to compute these in the
order of a few days on current technology. However, the algorithm is
highly parallelizable, and also, the intermediate output used
to estimate $|\Sigma(n,t)|$ can be used to ``bootstrap'' estimates of
$|\Sigma(n,u)|$ for $u>t$. Exploiting these features substantially
decrease computing time. Finally, the memory requirements of
our implementation are modest and have not been an issue in our experiments.

\subhead (1.5) Organization of the paper \endsubhead

 In Section 2, we describe the density $f$ of formula (1.3.1) 
and the factorization $f=p \phi$ of (1.3.2). 

In Section 3, we describe the algorithm of approximating the number 
$|\Sigma(n,t)|$ of magic squares in detail and state all the bounds for $|\Sigma(n,t)|$, $f$, $\phi$, and $p$ we need to conclude that the algorithm indeed approximates the desired number within relative error 
$\epsilon$ in $\left(1/\epsilon\right)^{O(1)} N^{O(\ln N)}$ time. We also 
describe details of the implementation. 

Sections 4-9 are devoted to the proofs.

 In Section 4, we invoke certain classical results, the van der Waerden and the Bregman-Minc bounds for the permanent of a non-negative matrix and obtain 
a straightforward corollary that we use later.

 In Section 5, we prove that the total number ${N +n^2-1 \choose n^2-1}$ 
 of $n \times n$ non-negative integer matrices with the sum
 of entries equal to $N$ is at most $N^{O(n)}$ times bigger than the number of $n \times n$ magic squares with the line sum $t$. Also, we prove that the maximum 
of the density $f$ on the simplex $\Delta$ in (1.3.1) does not exceed ${N +n^2-1 \choose n^2-1}$.

In Sections 6--8, we prove the key estimate of the paper, namely, that for any $\alpha>0$ there
is a $\beta=\beta(\alpha)>0$ such that the probability that a random $X \in \Delta$ satisfies
$p(X)> N^{\beta \ln N}$ in (1.3.2) does not exceed $N^{-\alpha n}$. Section 7 contains some standard probabilistic
estimates whereas Section 6 contains an estimate of the entries of the doubly stochastic scaling
of a positive matrix, which may be of interest in its own right. Roughly, it states that for a sufficiently
generic $n \times n$ matrix $A$, all the entries of its doubly stochastic scaling are sufficiently close to
$1/n$. 

In Section 9 we state some technical estimates for the log-concave density $\phi$ 
which imply that the algorithms of \cite{AK91}, \cite{F+94}, \cite{FK99}, \cite{LV06} for polynomial
time integration and sampling are indeed applicable. 

Finally, in Section 10 we describe possible
extensions of our approach, in particular, to contingency tables with equal row sums, but 
not necessarily column sums and vice versa. We also conjecture that the $N^{O(\ln N)}$ 
bound can be replaced by $N^{O(1)}$ so that our approach produces a polynomial time 
algorithm.

\head 2. The integral representation for the number of magic squares \endhead

To obtain the representation (1.3.1), we express the number $|\Sigma(n,t)|$ as the 
expectation of the permanent of an $N \times N$ random matrix for $N=nt$. 
Let $A=\left(a_{ij}\right)$ be an $N \times N$ 
square matrix. The {\it permanent} of $A$ is given by the formula
$$\per A=\sum_{\sigma \in S_N} \prod_{i=1}^N a_{i \sigma(i)},$$
where $S_N$ is the symmetric group of all permutations 
$\sigma$ of the set $\{1, \ldots, N\}$.  

Let $J_t$ denote the $t \times t$ matrix with the entries $1/t$. For an $n \times n$ matrix 
$X=\left(x_{ij}\right)$, the matrix $X \otimes J_t$ denotes the $N \times N$ block matrix 
whose $(i,j)$th block is the $t \times t$ matrix with the entries $x_{ij}/t$.

We recall that a random variable $\xi$ has the {\it standard exponential} distribution if 
$$\PP\bigl\{ \xi > \tau \bigr\}=\cases 1  &\text{if\ } \tau \leq 0 \\ e^{-\tau} &\text{if\ } \tau >0. \endcases$$ 

The following result is a particular case of the general formula of Theorem 1.2 of  \cite{Ba05}
 and Theorem 4 of \cite{Ba07} for the number of contingency tables of a given type.

\proclaim{(2.1) Theorem} Let $X=\left(x_{ij}\right)$ be the $n \times n$ matrix of independent 
standard exponential random variables $x_{ij}$. Then 
$$|\Sigma(n, t)| = {t^N \over (t!)^{2n}}\EE \per \left(X \otimes J_t \right).$$
\endproclaim

In other words, 
$$|\Sigma(n,t)|={t^N \over (t!)^{2n}} \int_{{\Bbb R}^{n^2}} \per \left(X \otimes J_t \right) 
\exp\left\{ - \sum_{i,j=1}^n x_{ij} \right\} \ d x, $$
where $dx$ is the standard Lebesgue measure in the space ${\Bbb R}^{n^2}$, interpreted 
as the space of $n \times n$ matrices $X=\left(x_{ij}\right)$. 

Let 
$$\Delta=\Delta_{n \times n} =\left\{ X=\left(x_{ij}\right): \quad x_{ij} >0 \quad 
\text{for all} \quad i,j \quad \text{and} \quad  \sum_{i,j=1}^n x_{ij}=1 \right\}$$
be the standard $(n^2-1)$-dimensional (open) simplex in ${\Bbb R}^{n^2}$ endowed with
the probability measure $\mu$ that is the normalization of the Lebesgue measure on 
$\Delta_{n \times n}$. 

Since $\per \left(X \otimes J_t \right)$ is a homogeneous polynomial of degree $N=nt$ in 
the entries $x_{ij}$ of $X$, we obtain the following integral representation, see Section 4 of \cite{Ba05}.
\proclaim{(2.2) Corollary}

$$|\Sigma(n, t)|= {(N+n^2 -1)! t^N \over (n^2-1)! (t!)^{2n}} \int_{\Delta_{n \times n}}
 \per \left(X \otimes J_t \right) \ d \mu(X).$$
\endproclaim

Thus we define $f: \Delta \longrightarrow {\Bbb R}_+$ in representation (1.3.1) by
$$f(X)={(N+n^2-1)!  t^N \over (n^2-1)! (t!)^{2n}} \per \left(X \otimes J_t \right) \quad \text{for}
\quad X \in \Delta_{n \times n}. \tag2.3$$

\subhead (2.4) Factoring the density \endsubhead
Now we describe how to factor the density
$f=p \phi$ in (1.3.2),
where $\phi: \Delta \longrightarrow {\Bbb R}_+$ is a log-concave function and 
$p: \Delta \longrightarrow {\Bbb R}_+$ is a function which ``does not vary much'' on 
$\Delta$. We employ the notion of {\it matrix scaling}, see \cite{Si64}, \cite{MO68}, \cite{KK96}, \cite{L+ 00}.

Let $X=\left(x_{ij}\right)$ be an $n \times n$ positive matrix. Then there exists an $n \times n$ positive
matrix $Y=\left(y_{ij}\right)$ and positive numbers $\lambda_i, \mu_i$ for $i=1, \ldots, n$ such that
$$x_{ij}=y_{ij} \lambda_i \mu_j \quad \text{for} \quad i,j=1, \ldots, n$$
and $Y$ is {\it doubly stochastic}, that is,
$$\split &\sum_{j=1}^n y_{ij}=1 \quad \text{for} \quad i=1, \ldots, n \quad \text{and} \\
&\sum_{i=1}^n y_{ij}=1 \quad \text{for} \quad j=1, \ldots, n. \endsplit$$

Furthermore, given matrix $X$, the matrix $Y$ is unique 
(we call it the {\it doubly stochastic scaling of} $X$) while the factors $\lambda_i$ and $\mu_j$ 
are unique up to a rescaling $\lambda_i :=\lambda_i \tau$, $\mu_i:=\mu_i \tau^{-1}$ 
for some $\tau>0$ and $i=1, \ldots, n$. This allows us to define a function 
$\sigma: \Delta \longrightarrow {\Bbb R}_+$ by 
$$\sigma(X)=\prod_{i=1}^n \left(\lambda_i \mu_i\right).$$
Clearly,
$\per X = \left(\per Y \right) \sigma(X)$.

The crucial fact about $\sigma$ that we use is that $\sigma$ is log-concave, that is,
$$\sigma\left(\alpha_1 X_1 +\alpha_2 X_2 \right) 
\geq \sigma^{\alpha_1}(X_1) \sigma^{\alpha_2}(X_2)$$
for all $X_1, X_2 \in \Delta$ and all $\alpha_1, \alpha_2 \geq 0$ such that $\alpha_1+\alpha_2=1$,
see \cite{GS02}, \cite{Gu06}, \cite{Ba05}, \cite{Ba06}. Also, for any given $X$, the value of $\sigma$ can be computed within relative error $\epsilon$ 
in time polynomial in $\ln  \left(1/\epsilon\right)$ and $n$ \cite{L+00}. 

One can easily see that if $Y$ is the doubly stochastic scaling of $X$ then $Y \otimes J_t$ is 
the doubly stochastic scaling of $X \otimes J_t$ and that
$$\sigma\left(X \otimes J_t \right)=\sigma^t(X).$$
Moreover,
$$\per \left(X \otimes J_t \right)=\per \left(Y \otimes J_t \right) \sigma^t(X).$$
We define 
$$\aligned  & f=p \phi \quad \text{where} \\ &p(X)={N^N \over N!} \per \left(Y \otimes J_t \right) \quad 
\text{and} \\ & 
\phi(X)= {(N+n^2-1)! N! t^N \over  (n^2-1)! (t!)^{2n} N^N}\sigma^t(X),
\endaligned \tag2.4.1$$
for $N=nt$. Clearly, $\phi$ is log-concave and Corollary 2.2 implies that
$$|\Sigma(n,t)|=\int_{\Delta} p \phi \ d \mu. \tag2.4.2$$
We note that for any given $X \in \Delta$ and $\epsilon>0$ the value of $p(X)$ can be computed
within relative error $\epsilon$ in time polynomial in $N$ and $\epsilon^{-1}$ by the randomized algorithm 
of \cite{J+04}.

It is convenient to define $f$ and $\phi$ on the set $\Mat_+$ of positive $n \times n$ matrices
and not just on the simplex $\Delta$.

\head 3. Bounds and the detailed description of the algorithm \endhead

In this section we give the detailed description of the algorithm and also summarize various bounds that we need. We begin with some general bounds
on the number $|\Sigma(n,t)|$ of magic squares and the density factors $p$ and $\phi$ in (1.3.2)
and (2.4.1).
\proclaim{(3.1) Theorem}
We have 
\roster
\item 
$${N+n^2-1 \choose n^2-1} \geq |\Sigma(n,t)| \geq 
{N+n-1 \choose n-1}^{-2} {N+n^2-1 \choose n^2-1};$$
\item
$$1 \leq p(X) \leq {(t! )^n N^N \over t^{N} N!} \quad \text{for all} \quad X \in \Delta;$$ 
\item
$$0 < \phi(X) \leq {(N+n^2-1)! N!  \over (n^2-1)! (t!)^{2n} n^{2N}}\quad \text{for all} \quad X \in \Delta.$$
\endroster
\endproclaim

We will use somewhat cruder estimates which are easier to work with.

\proclaim{(3.2) Corollary} We have
\roster
\item
$${N+n^2-1 \choose n^2-1} \geq  |\Sigma(n,t)| \geq \left(N+n\right)^{-2n}{N+n^2-1 \choose n^2-1};$$
\item
$$ 1 \leq p(X) \leq (8 t)^{n/2} \quad \text{for all} \quad X \in \Delta;$$
\item
$$ 0 \leq f(X) \leq {N+n^2-1 \choose n^2-1} \quad \text{for all} \quad X \in \Delta.$$
\endroster
\endproclaim

We prove Theorem 3.1 in Section 5. Corollary 3.2 follows by standard estimates via Stirling's formula
$$\sqrt{2 \pi x} \left({x \over e}\right)^x e^{1/(12x +1)} \leq \Gamma(x+1) \leq 
\sqrt{2 \pi x} \left({x \over e}\right)^x e^{1/(12 x)} \quad \text{for} \quad x \geq 0. \tag3.2.1$$

It turns out that most of the time we have $p(X)=N^{O(\ln N)}$.
Now we state the key estimate of the paper.

\proclaim{(3.3) Theorem} For any $\alpha>0$ there exists $\beta=\beta(\alpha)>0$ such that for all positive integers 
$n$ and $t$ such that 
$$t < e^n,$$
we have
$$\mu\Bigl\{ X \in \Delta_{n \times n}: \quad p(X) \geq N^{\beta \ln N}  \Bigr\} < N^{-\alpha n}.$$
\endproclaim

We prove Theorem 3.3 in Section 8 having established an estimate for the entries of the doubly 
stochastic scaling of a matrix in Section 6 and some standard probability bounds in Section 7. 

Finally, we need some technical estimates showing that $\phi$ is sufficiently regular so that we can indeed apply integration and sampling algorithms of \cite{AK91},
\cite{F+94}, \cite{FK99}, and \cite{LV06}.

For 
$$0 < \delta < {1 \over n^2},$$
let us define the {\it $\delta$-interior} of the simplex by 
$$\Delta^{\delta}=\Delta_{n \times n}^{\delta}=\left\{X=\left(x_{ij}\right): \quad
x_{ij}>\delta \quad \text{for all} \quad i,j \quad \text{and} \quad \sum_{i,j=1}^n x_{ij}=1 \right\}.$$

\proclaim{(3.4) Theorem} We have
\roster
\item
$$\int_{\Delta^{\delta}} f \ d \mu \leq \int_{\Delta} f \ d \mu \leq \left(1- \delta n^2\right)^{-N-n^2+1} 
\int_{\Delta^{\delta}} f \ d \mu,  $$
\item 
For any $X, Y \in \Delta^{\delta}$ where $X=\left(x_{ij}\right)$ and $Y=\left(y_{ij}\right)$ we have 
$$| \ln \phi(X) -\ln \phi(Y) | \leq {N \over \delta} \max_{ij} |x_{ij}-y_{ij}|.$$
\endroster
\endproclaim

Theorem 3.4 is proven in Section 4 of \cite{Ba05}. For completeness, we present its straightforward proof
in Section 9.

\subhead (3.5) The algorithm \endsubhead Now we can describe the algorithm in more detail.
First, we assume that $t< e^n$. Indeed, for $t \geq n^3$ there is a randomized polynomial time approximation scheme for computing $|\Sigma(n,t)|$ which is a particular case of the algorithm of 
Dyer, Kannan, and Mount  \cite{D+97}, see also \cite{Mo02} for a strengthening.
Second, we assume that $\epsilon > N^{-n}$, which is not really restrictive since $|\Sigma(n, t)|$ 
can be computed exactly in $N^{O(n)}$ time by a dynamic programming algorithm.

For an $0< \epsilon<1$, let us choose
$$\delta ={-\ln(1-\epsilon) \over n^2(N+n^2-1)} \approx {\epsilon \over n^2(N+n^2-1)} \quad
\text{for small} \quad \epsilon>0.$$
By Part 1 of Theorem 3.4, the integral
$$\int_{\Delta^{\delta}} f \ d\mu \tag3.5.1$$
approximates $|\Sigma(n,t)|$ from below within the relative error $\epsilon$. We factor 
$f=p\phi$ as in Section 2.4. By Parts 1 and 2 of Corollary 3.2 and Theorem 3.3 (with a sufficiently large 
$\alpha$)
it follows that for $\overline{p}$ defined by
$$\overline{p}(X)=\cases p(X) &\text{if} \quad p(X) \leq T \\ T &\text{if} \quad p(X)>T
\endcases$$
with $T=N^{\beta \ln N}$,
the integral 
$$\int_{\Delta^{\delta}}  \overline{p} \phi \ d \mu \tag3.5.2$$
approximates (3.5.1) from below within relative error $N^{-n} < \epsilon$.

A simplified, polynomial time algorithm, replaces integral (3.5.2) by the integral
$$\int_{\Delta^{\delta}} \phi \ d \mu, \tag3.5.3$$
which approximates (3.5.2) within a factor of $N^{\beta \ln N}$. Because $\phi$ is log-concave
and satisfies the bound of Part 2 of Theorem 3.4, the integral (3.5.3) can be computed within relative 
error $\epsilon$ by any of the algorithms of \cite{AK91},
\cite{F+94}, \cite{FK99}, and \cite{LV06} in time polynomial in $N$ and $\epsilon^{-1}$. We use 
an algorithm of \cite{L+00} to compute $\phi(x)$ for a given $x \in \Delta$.

In the more accurate, but also more time consuming, version of the algorithm, 
we write (3.5.2) as 
$$\left(\int_{\Delta^{\delta}} \overline{p} \ d \nu\right) \left(\int_{\Delta^{\delta}} \phi \ d \mu \right),$$
where $\nu$ is the probability measure on $\Delta^{\delta}$ with the density proportional to
$\phi$. 
Furthermore, the algorithms of \cite{AK91},
\cite{F+94}, \cite{FK99}, and \cite{LV06} allow us to sample independent random 
points from a probability measure $\tilde{\nu}$ sufficiently close to $\tilde{\nu}$, that is 
satisfying 
$$\left|\tilde{\nu}(S)- \nu(S) \right| < \epsilon N^{-\beta \ln N} \quad \text{for any Borel}
\quad S \subset \Delta^{\delta}.$$
A single point can be sampled in 
$ \left(1/\epsilon\right)^{O(1)} N^{O(\ln N)}$ time. We sample 
$m=\lceil 3T^2 \epsilon^{-2} \rceil$ independent random points $x_i$ with respect to measure $\tilde{\nu}$
and estimate 
the integral 
$$\int_{\Delta^{\delta}} \overline{p} \ d \tilde{\nu} \tag3.5.4$$
by the sample mean
$$m^{-1} \sum_{i=1}^m \overline{p}\left(x_i\right) \tag3.5.5$$
By Chebyshev inequality, (3.5.5) approximates (3.5.4) within relative error $\epsilon$ 
with probability at least $2/3$. We use the algorithm \cite{J+04} to compute $p(x_i)$.

\subhead (3.6) Details of the implementation \endsubhead We implemented a much simplified version 
of the algorithm, computing the integral (1.3.4)
$$\int_{\Delta} \phi \ d \mu.$$
In our implementation, we work with the original simplex $\Delta$, not its $\delta$-interior 
$\Delta^{\delta}$. This has never given us any boundary-related trouble in our computational
experiments.

The implementation is based on a version of the hit-and-run algorithm
of \cite{LV06}, see also \cite{Ve05}. We use telescoping with respect to the density $\phi$.
Namely, we pick a sufficiently dense uniform subset
$$0 < t_1 < t_2 < \ldots < t_m =t,$$
and define a log-concave function $\psi_i$ by 
$$\psi_i(X)=\sigma^{t_i}(X) \quad \text{for} \quad X \in \Delta,$$
cf. Section 2.4 and formula (2.4.1), in particular. The number $m$ of points is chosen by the 
user so as to be ``reasonable''.

For a given $X \in \Delta$ we compute $\sigma(X)$ by the Sinkhorn balancing (alternate 
scaling of rows and columns of $X$ to the unit sums) \cite{Si64}, 
which seems to work very well even for $n$ as large as $n=100$.

 Note that the numbers $t_i$ are not necessarily integer 
and that 
$$\phi(X)={\left(N+n^2-1\right)! N! t^{N} \over (n^2-1)! \left(t!\right)^{2n} N^N} \psi_m(X),$$
cf. (2.4.1).
Hence the goal is to compute 
$$\split \int_{\Delta} \psi_m \ d \mu=S_1 \prod_{i=1}^{m-1}
&{S_{i+1} \over S_i} \quad \text{where} \\
&S_i=\int_{\Delta} \psi_i \ d \mu \quad \text{for} \quad i=1, \ldots, m. \endsplit$$

If $t_1$ is sufficiently small, the function $\psi_1$ is close to a constant, so we estimate 
$S_1$ by a sample mean of $\psi_1$ for a set of randomly chosen $X \in \Delta$.  In our experiments,
we often chose $t_1=1$. To choose 
a random $X \in \Delta$ from $\mu$, we choose the entries $x_{ij}$ of $X$ independently from the 
standard exponential distribution and then normalize:
$$x_{ij}:=x_{ij} \left( \sum_{k,l=1}^n x_{kl} \right)^{-1} \quad \text{for} \quad i,j=1, \ldots, n.$$

To compute ratios 
$${S_{i+1} \over S_i}=\int_{\Delta} {\psi_{i+1} \over \psi_i} \ d \nu_i,$$
where $\nu_i$ is the probability measure on $\Delta$ with the density proportional to $\psi_i$,
 we sample points $X \in \Delta$ from $\nu_i$ and average the ratios  $\psi_{i+1}(X)/\psi_i(X)$. If $t_{i+1}$ is sufficiently 
close to $t_i$, the ratios are close to 1, so a moderate number of sampled points  is needed.
Again, the number of samples is selected by the user.
Thus the bottleneck of this simplified
 algorithm consists in sampling random points $X \in \Delta$ from the 
probability measure $\nu_i$. For that, we iterate the basic ``hit-and-run'' construction.
We sample a random point $X_0 \in \Delta$ from the uniform distribution $\mu$, pick a random 
line $\ell$ through $X$ in the affine hull of $\Delta$, sample $X_1$ from the density on the 
interval $\ell \cap \Delta$ proportional to the restriction of $\psi_i$ onto that interval and 
iterate the process with $X_0:=X_1$ as the new starting point. After a number of iterations chosen 
by the user, the point $X_1$ is accepted as a random sample from the measure $\nu_i$.

To choose a random line $\ell$ through $X$ we first sample an $n \times n$ matrix 
$L=\left(\lambda_{ij}\right)$ of independent standard Gaussian random variables and then 
center it
$$\lambda_{ij}:=\lambda_{ij}-{1 \over n^2} \sum_{k,l=1}^n \lambda_{kl} \quad \text{for} \quad i,j=1, \ldots, n.$$
We then define the line $\ell$ through $X$ by 
$$\ell=\Bigl\{X + \tau L: \quad \tau \in {\Bbb R} \Bigr\}.$$
To choose a point $X_1 \in \ell \cap \Delta$, we approximate the restriction of $\psi_i$ onto
$\ell \cap \Delta$ by a function $\widetilde{\psi_{i, \ell}}$ such that $\ln \widetilde{\psi_{i, \ell}}$ is 
a piece-wise linear approximation of the restriction of $\ln \psi_i$ onto $\ell \cap \Delta$.
Sampling from the density proportional to $\widetilde{\psi_{i, \ell}}$ reduces then to sampling 
from the exponential distribution.

\head 4. Preliminaries: estimates on the permanent \endhead

We will use the following bounds for the permanent.

\subhead (4.1) The van der Waerden bound \endsubhead Let $B=\left(b_{ij}\right)$ be an 
$N \times N$ doubly stochastic matrix, that is, 
$$\sum_{j=1}^N b_{ij}=1 \quad \text{for} \quad i=1, \ldots, N \quad \text{and} 
\quad \sum_{i=1}^N b_{ij}=1 \quad \text{for} \quad j=1, \ldots, N$$
and 
$$b_{ij} \geq 0 \quad \text{for} \quad i,j=1, \ldots, N.$$
Then
$$ \per B \geq {N! \over N^N}.$$
This is the famous van der Waerden bound proved by Falikman \cite{Fa81} and Egorychev
\cite{Eg81}, see also Chapter 12 of \cite{LW01}.

\subhead (4.2) The continuous version of the Bregman-Minc bound \endsubhead
Let $B=\left(b_{ij}\right)$ be an $N \times N$ matrix such that 
$$\sum_{j=1}^N b_{ij} \leq 1 \quad \text{for} \quad i=1, \ldots, N$$
and 
$$b_{ij} \geq 0 \quad i,j=1, \ldots, N.$$
Furthermore, let 
$$s_i =\max_{j=1, \ldots, N} b_{ij}>0 \quad \text{for} \quad i=1, \ldots, N.$$
Then 
$$\per B \leq \prod_{i=1}^N s_i \Gamma^{s_i} \left({1+s_i \over s_i}\right).$$
This bound was obtained by Soules \cite{So03}.

If $s_i=1/r_i$ for  integers $r_i$ the bound transforms into 
$$\per B \leq \prod_{i=1}^N {(r_i!)^{1/r_i} \over r_i},$$
which can be easily deduced from the Minc conjecture proved by  Bregman 
\cite{Br73}, see also Chapter 11 of \cite{LW01}.

We will use the following corollary of estimates of Sections 4.1 -- 4.2.

\proclaim{(4.3) Corollary} Let $B=\left(b_{ij}\right)$ be an $N \times N$ doubly stochastic 
matrix and let 
$$s_i=\max_{j=1, \ldots, N} b_{ij} \quad \text{for} \quad i=1, \ldots, N.$$
Suppose that
$$\sum_{i=1}^N s_i \leq \gamma \quad \text{for some} \quad \gamma \geq 1.$$
Then 
$$\split {N! \over N^N} \leq \per B \leq &\left({ \gamma \over N} \right)^N \Gamma^{\gamma} \left(1+{N \over \gamma} \right) \\ \leq & {N! \over N^N} \left(2 \pi N \right)^{\gamma/2} e^{\gamma^2/12 N} . 
\endsplit.$$
\endproclaim
\demo{Proof}
The lower bound is the van der Waerden estimate, see Section 4.1.

Let us define
$$g(\xi)=\xi\ln \Gamma \left({1+\xi \over \xi}\right)  +\ln \xi \quad \text{for} \quad 0< \xi \leq 1.$$
Then $g$ is a concave function, cf.  \cite{So03}, and by the inequality of Section 4.2, we have 
$$\ln \per B \leq \sum_{i=1}^N g(s_i).$$
The function 
$$G(x) =\sum_{i=1}^N g(\xi_i) \quad \text{for} \quad x=\left(\xi_1, \ldots, \xi_N\right)$$ 
is concave on the simplex defined by the equation
$\xi_1 + \ldots + \xi_N = \gamma$ and inequalities $\xi_i \geq 0$ for $i=1, \ldots, N$.
It is also symmetric under permutations of $\xi_1, \ldots, \xi_N$. 
Hence the maximum of $G$ is attained at 
$$\xi_1 = \ldots = \xi_N =\gamma/N ,$$
and so 
$$\ln \per B \leq N g\left({\gamma \over N} \right).$$
Thus
$$\per B \leq  \left({ \gamma \over N} \right)^N \Gamma^{\gamma} \left(1+{N \over \gamma} \right)$$
and the rest follows by Stirling's formula (3.2.1).
{\hfill \hfill \hfill} \qed
\enddemo
We will apply Corollary 4.3 for $\gamma=O(\ln N)$, in which case the ratio of the
upper and lower bounds is $N^{O(\ln N)}$.

\head 5. Proof of Theorem 3.1 \endhead

To prove the upper bound in Part  1, we note that $|\Sigma(n,t)|$ does not exceed the number 
of $n \times n$ non-negative integer matrices with the sum of entries equal to $N$, which is 
exactly equal to ${N +n^2-1 \choose n^2-1}$.

For non-negative integer vectors $R=\left(r_1, \ldots, r_n\right)$ and $C=\left(c_1, \ldots, c_n\right)$
such that $r_1 + \ldots + r_n =c_1 + \ldots + c_n=N$, let $|\Sigma(R, C)|$ be the number of 
$n \times n$ non-negative integer matrices with the row sums $r_1, \ldots, r_n$ and the
column sums $c_1, \ldots, c_n$, cf. Section 1.2.
 As is discussed in \cite{Ba07}, see formula (4) of Section 2 there, we have
$$|\Sigma(R,C)| \leq |\Sigma(n,t)| \quad \text{for all} \quad R \quad \text{and} \quad C.$$
Since the the total number of pairs $(R, C)$ is ${N + n-1 \choose n-1}^2$, the lower bound 
follows by summing up the previous inequality over all choices of $(R,C)$.

We recall that  $p$ is defined by
$$p(X)={N^N \over N!} \per \left(Y \otimes J_t \right),$$
where $Y$ is an $n \times n$ doubly stochastic scaling of $X$ and $J_t$ is the $t \times t$ matrix filled 
with $1/t$, cf. Section 2.4. Hence $Y \otimes J_t$ is an $N \times N$ doubly stochastic matrix and the lower 
bound in Part 2 follows by the van der Waerden bound, see Section  4.1. Since the entries of 
$Y \otimes J_t$ do not exceed $1/t$, by Section 4.2, we have
$$p(x) \leq {N^N \over N!} {(t!)^n \over t^N},$$
which completes the proof of Part 2.

Clearly, 
$$\phi(X) ={(N+n^2-1)! N! t^N \over (n^2-1)!(t!)^{2n} N^N} \sigma^t(X) > 0 \quad \text{for all} 
\quad X \in \Delta.$$
Since $\sigma(X)$ is log-concave and invariant under permutations of rows and columns of 
$X$, the maximum of $\sigma$ on $\Delta$ is attained at the $n \times n$ matrix $A$ with the entries 
$1/n^2$ that is the average of all the matrices obtained from any given $X \in \Delta$ by permutations
of rows and columns. It is immediate that $\sigma(A)=n^{-n}$ and hence 
$$\phi(X) \leq \phi(A)={(N+n^2-1)! N! t^N \over (n^2-1)! (t!)^{2n} N^N} n^{-N},$$
and the proof of Part 3 follows.
{\hfill \hfill \hfill} \qed

As we remarked before, Corollary 3.2 is obtained by a straightforward application of the 
Stirling formula and a remark that $f=p\phi$.

\head 6. Bounding the entries of the doubly stochastic scaling of a matrix \endhead

In this section we prove the following main result.

\proclaim{(6.1) Theorem} Let $A=\left(a_{ij}\right)$ be an $n \times n$ positive matrix 
and let $B=\left(b_{ij}\right)$  be the doubly stochastic scaling of $A$
so that for some $\lambda_i, \mu_j >0$ we have 
$$ \split &a_{ij}=b_{ij} \lambda_i \mu_j \quad \text{for all} \quad i,j \quad \text{and} \\
&\sum_{j=1}^n b_{ij}=1 \quad \text{for} \quad i=1, \ldots, n \quad \text{and} \quad 
\sum_{i=1}^n b_{ij}=1 \quad \text{for} \quad j=1, \ldots, n. \endsplit$$
Then, for all $k$ and $l$,
$$\split  \ln b_{kl} \leq \ln a_{kl} &-{1 \over n-2} \sum_{j \ne l} \ln a_{kj}-{1 \over n-2} \sum_{i\ne k}
\ln a_{il}\\&+{n \over n-2}\ln\left({1 \over n} \sum_{i,j=1}^n
 a_{ij} \right) -{2n-2 \over n-2}\ln (n-1).\endsplit $$
\endproclaim

\example{Example} Suppose that $1 \leq a_{ij} \leq 2$ for all $i,j$.
Theorem 6.1 implies that for some absolute constant $\gamma$ we have
$$b_{ij} \leq {\gamma \over n} \quad \text{for all} \quad i,j.$$
Hence by Corollary 4.3, 
$${n! \over n^n} \leq \per B \leq n^{O(1)} {n! \over n^n}$$
and hence values of $\per B$ vary within up to a polynomial in $n$ factor. In contrast,
$$n! \leq \per A \leq 2^n n!,$$
so values of $\per A$ vary within an exponential in $n$ factor.

This concentration of the permanent of the doubly stochastic scaling of a matrix is the basis 
of our approach.
\endexample

The proof of Theorem 6.1 is based on the following two lemmas.

The first lemma was proved in \cite{L+00}, for completeness we include its  proof here.
\proclaim{(6.2) Lemma} Let $A=\left(a_{ij}\right)$ be an $n \times n$ positive matrix 
such that 
$$\sum_{i,j=1}^n  a_{ij}=n$$ and let $B=\left(b_{ij}\right)$ be the doubly stochastic scaling of $A$.
Then
$$\sum_{i,j=1}^n \ln b_{ij} \geq \sum_{i,j=1}^n \ln a_{ij}.$$
\endproclaim
\demo{Proof} As is known, $B$ can be obtained from $A$ as the limit of repeated 
alternate scalings of the rows of $A$ to the row sums equal to 1 and of the columns of $A$ to the column sums equal to 1 (Sinkhorn balancing), see \cite{Si64}.
Hence it suffices to prove that under the row (column) scalings, the sum of the logarithms 
of the entries of the matrix can only increase.

To this end, let 
$C=\left(c_{ij}\right)$ be a positive $n \times n$ matrix with the row sums 
$\rho_1, \ldots, \rho_n$ such that $\rho_1 + \ldots + \rho_n=n$ and 
let $D=\left(d_{ij}\right)$ be the matrix such that 
$$d_{ij}={c_{ij} \over \rho_i} \quad \text{for all} \quad i,j.$$
In words: we divide the $i$th row of $C$ by its row sum $\rho_i$. We note that the sum of the 
entries of $D$ is $n$ and that 
$$\sum_{i,j=1}^n  \left(\ln d_{ij} - \ln c_{ij}\right) =-n\sum_{i=1}^n \ln \rho_i \geq 0$$
because of the arithmetic-geometric mean inequality.

Column scalings are handled in the same way.
{\hfill \hfill \hfill} \qed
\enddemo

The following result is obtained in \cite{Br73}. For completeness, we 
provide its proof below.

\proclaim{(6.3) Lemma} Let $A=\left(a_{ij}\right)$ be a positive $n \times n$ matrix and 
let $B=\left(b_{ij}\right)$ be its doubly stochastic scaling. Then $B$ is the solution of the 
optimization problem 
$$\text{minimize} \quad \sum_{i,j=1}^n x_{ij} \left( \ln x_{ij} - \ln a_{ij}\right)$$
over the set of all $n \times n$ doubly stochastic matrices $X=\left(x_{ij}\right)$.
\endproclaim
\demo{Proof} First, we note the minimum is attained on a positive doubly stochastic
matrix $X$. If, for example, $x_{11}=0$, then there are indices $i$ and $j$ such that 
$x_{i1}>0$, $x_{1j}>0$, $x_{ij}<1$ and one can make the value of the objective function 
smaller by modifying
$$x_{11}:=\epsilon, \quad x_{1j}:=x_{1j}-\epsilon, \quad x_{i1}:=x_{i1}-\epsilon, \quad
x_{ij}:=x_{ij}+\epsilon$$
for a sufficiently small $\epsilon >0$. This follows since the right derivative of $x \ln x$ 
is $-\infty$ at $x=0$ and is finite at any $x>0$.

Since the optimal point $X$ lies in the relative interior of the set of doubly stochastic 
matrices, the gradient of the objective function at $X$ should be orthogonal to the space of 
$n \times n$ matrices with the row and column sums equal to $0$.

This gives us the following equations
$$\ln x_{ij}-\ln a_{ij} =\xi_i + \eta_j\quad \text{for all} \quad i,j$$
and some numbers $\xi_1, \ldots, \xi_n; \eta_1, \ldots, \eta_n$.

In other words,
$$x_{ij} =a_{ij} \lambda_i \mu_j \quad \text{for} \quad \lambda_i=e^{\xi_i}
\quad \text{and} \quad \mu_j= e^{\eta_j} \quad \text{for all} \quad i,j$$
as desired.
{\hfill \hfill \hfill} \qed
\enddemo

Now we are ready to prove Theorem 6.1.

\demo{Proof of Theorem 6.1}  First, we notice that  neither the matrix $B$ nor the right hand side of the 
inequality change if we scale 
$$a_{ij}:=a_{ij} \tau \quad \text{for all}\quad i,j$$
and some $\tau>0$. Therefore, without loss of generality, we assume that 
$$\sum_{i,j=1}^n a_{ij} =n.$$

Without loss of generality, we assume that $k=l=1$, so our 
goal is to bound $b_{11}$.

By Lemma 6.3, matrix $B$ is the solution of the minimization problem
$$\text{minimize} \quad \sum_{i,j=1}^n x_{ij}\left(\ln x_{ij}- \ln a_{ij}\right)$$
over the set of $n \times n$ doubly stochastic matrices $X=\left(x_{ij}\right)$.

For a real $\tau$, let us define the matrix $B(\tau)=\left(b_{ij}(\tau)\right)$ by
$$b_{ij}(\tau)=\cases b_{11}+\tau  &\text{if\ } i=j=1, \\
b_{1j}-\tau/(n-1) &\text{if\ }i=1, j \ne 1,\\
b_{i1}-\tau/(n-1) &\text{if \ } i \ne 1, j=1, \\
b_{ij}+\tau/(n-1)^2 &\text{if \ } i \ne 1, j \ne 1. \endcases$$

We observe that $B(0)=B$, that the row and column sums of $B(\tau)$ are 1 and that $B(\tau)$ is 
positive for all $\tau$ 
from a sufficiently small neighborhood of the origin.
Therefore, if we let 
$$f(\tau)=\sum_{i,j=1}^n b_{ij}(\tau)\left(\ln b_{ij}(\tau)-\ln a_{ij}\right),$$
we must have 
$$f'(0)=0.$$

Computing the derivative, we get 
$$\split f'(0)=&\ln b_{11}- \ln a_{11} + 1 \\
                       &-{1 \over n-1} \sum_{j \ne 1} \left(\ln b_{1j}-\ln a_{1j}\right)-1 \\
                        &-{1 \over n-1} \sum_{i \ne 1} \left(\ln b_{i1}- \ln a_{i1} \right)  -1\\  
                        &+{1 \over (n-1)^2} \sum_{i,j \ne 1}  \left(\ln b_{ij} -\ln a_{ij}\right) +1.
  \endsplit$$       
Rearranging summands, we rewrite the derivative in the form
$$\split f'(0)=&\left(1-{1 \over (n-1)^2} \right) \left(\ln b_{11}-\ln a_{11} \right) \\
                      &-\left({1 \over n-1} +{1 \over (n-1)^2} \right) \sum_{j \ne 1} \left(\ln b_{1j}-\ln a_{1j}\right)\\
                     & -\left({1 \over n-1} +{1 \over (n-1)^2} \right) \sum_{i \ne 1} \left(\ln b_{i1}-\ln a_{i1}\right)\\
                      &+{1 \over (n-1)^2} \sum_{i,j=1}^n  \left(\ln b_{ij} -\ln a_{ij}\right). \endsplit$$
Since $f'(0)=0$ and since  by Lemma 6.2 we have                      
$$\sum_{i,j=1}^n \ln b_{ij} \geq \sum_{i,j=1}^n \ln a_{ij},$$
we must have 
$$\split &{n^2-2n \over (n-1)^2}  \left(\ln b_{11}-\ln a_{11} \right)  \\ 
              &- {n \over (n-1)^2} \sum_{j \ne 1} \left(\ln b_{1j}-\ln a_{1j}\right) 
              -{n \over (n-1)^2}  \sum_{i \ne 1} \left(\ln b_{i1}-\ln a_{i1}\right) \\
              & \leq 0.\endsplit$$  
 That is,             
  $$(n-2)   \left(\ln b_{11}-\ln a_{11} \right)  - \sum_{j \ne 1} \left(\ln b_{1j}-\ln a_{1j}\right) 
              - \sum_{i \ne 1} \left(\ln b_{i1}-\ln a_{i1}\right) \leq 0.$$             
 In other words,
 $$\split & \ln b_{11}  -   {1 \over n-2} \sum_{j \ne 1} \ln b_{1j}   -{1 \over n-2}
 \sum_{i \ne 1} \ln b_{i1} \\ \leq &\ln a_{11}  -   {1 \over n-2} \sum_{j \ne 1} \ln a_{1j}   -{1 \over n-2}
 \sum_{i \ne 1} \ln a_{i1} \endsplit$$
 and
 $$\split  \ln b_{11}  \leq &\ln a_{11}  -   {1 \over n-2} \sum_{j \ne 1} \ln a_{1j}   -{1 \over n-2}
 \sum_{i \ne 1} \ln a_{i1} \\  &+     {1 \over n-2} \sum_{j \ne 1} \ln b_{1j}   +{1 \over n-2}
 \sum_{i \ne 1} \ln b_{i1}. \endsplit$$
 On the other hand, if the value of $b_{11}$ is fixed, the maximum value of 
 $$\sum_{j \ne 1} \ln b_{1j} +\sum_{i \ne 1} \ln b_{i1} $$
 is attained at 
 $$b_{1j}=b_{i1}={1-b_{11} \over n-1} \quad \text{for all} \quad i,j \ne 1$$
 (since the row and column sums of $B$ are equal to 1).
 
 Therefore, we have
 $$\split \ln b_{11} \leq  &\ln a_{11}  -   {1 \over n-2} \sum_{j \ne 1} \ln a_{1j}   -{1 \over n-2}
 \sum_{i \ne 1} \ln a_{i1} \\ 
  & +{2n-2 \over n-2} \ln\left(1-b_{11}\right) -{2n-2 \over n-2}\ln(n-1).\endsplit$$
  
 Since $\ln\left(1-b_{11}\right) \leq 0$ and $\sum_{i,j=1}^n a_{ij}=n$, this completes the proof.
 {\hfill \hfill \hfill \qed}
\enddemo

\head 7. Probabilistic estimates \endhead

The goal of this section is to prove the following technical estimates.

\proclaim{(7.1) Lemma} 

\roster

\item Let $a_1, \ldots, a_n$ be independent standard exponential random variables. 
Then for every $r>0$ we have 
$$\PP\left\{ {1 \over n} \sum_{i=1}^n \ln a_i \leq -r \right\} \leq \left({\pi \over e^r}\right)^{n/2};$$

\item Let $a_1, \ldots, a_m$ be independent standard exponential random variables.
Then for every $r>0$ we have 
$$\PP\left\{ {1 \over m} \sum_{i=1}^m a_i \geq r \right\} \leq \left({4 \over e^r}\right)^{m/2}.$$

\item Let $a_{ij}$, $1 \leq i,j \leq n$, be independent standard exponential random variables
and let 
$$c_i=\max_{j=1, \ldots, n} a_{ij} \quad \text{for} \quad i=1, \ldots, n.$$
Then for every $r>0$ we have
$$\PP\left\{ {1 \over n} \sum_{i=1}^n c_i > r \right\} \leq e^{-rn/2} \left(n \sqrt{\pi}\right)^n.$$
\endroster

\endproclaim
\demo{Proof} We use the Laplace transform method, see, for example, Appendix A of  \cite{AS00}.

To prove Part 1,
let
$$b={1 \over n} \sum_{i=1}^n \ln a_i.$$
For $\tau>0$ we get
$$\PP\left\{b < -r \right\}=\PP\left\{e^{-\tau b} > e^{\tau r} \right\} \leq e^{-\tau r} \EE e^{-\tau b}$$ 
by the Markov inequality.
Let us choose $\tau=n/2$. Then 
$$\EE e^{-\tau b}=\EE \prod_{i=1}^n a_i^{-1/2}=\left( \int_0^{+\infty} x^{-1/2} e^{-x} \ dx \right)^n
=\Gamma^n(1/2)=\pi^{n/2}.$$

Hence
$$\PP\left\{ {1 \over n} \sum_{i=1}^n \ln a_i < -r \right\} \leq \pi^{n/2} e^{-rn/2}.$$

To prove Part 2, let 
$$b={1 \over m}\sum_{i=1}^m a_i.$$
For $\tau>0$ we get
$$\PP\left\{ b>r\right\} =\PP\left\{ e^{\tau b} > e^{\tau r} \right\} \leq e^{-\tau r} \EE e^{\tau b}.$$
Let us choose $\tau=m/2$.
Then
$$\EE e^{\tau b} =\prod_{i=1}^m \EE e^{a_i/2}= \left( \int_0^{+\infty} e^{-x/2} \ dx \right)^m=2^m.$$

To prove Part 3, let
$$b={1 \over n} \sum_{i=1}^n c_i.$$ 
For $\tau>0$ we get
$$\PP\left\{ b>r \right\} \leq e^{-\tau r} \EE e^{\tau b}.$$
We choose $\tau=m/2$. Then 
$$\EE e^{\tau b}=\prod_{i=1}^n \EE e^{c_i/2}.$$

We have
$$\split &\EE \left(e^{c_i/2}\right) = \EE \left(\max_{j=1, \ldots, n} e^{a_{ij}/2} \right) \leq 
\EE \left( \sum_{j=1}^n e^{a_{ij}/2} \right) \\
  =&n \int_0^{+\infty}  e^{-x/2} \ dx =n \Gamma(1/2) =n \sqrt{\pi}\endsplit $$
and the proof follows.
{\hfill \hfill \hfill} \qed
\enddemo

\head 8. Proof of Theorem 3.3 \endhead

Let $\Mat_+=\Mat_+(n,n)$ be the set of $n \times n$ positive matrices $A=\left(a_{ij}\right)$ and let us consider the projection $\Psi: \Mat_+(n, n) \longrightarrow \Delta_{n \times n}$, where 
$\Psi(A)=X=\left(x_{ij}\right)$ is defined by
$$x_{ij}=a_{ij}\left(\sum_{k,l=1}^n a_{kl}\right)^{-1} \quad \text{for} \quad i,j=1, \ldots, n.$$
As is known and easy to check, the push-forward of 
the exponential measure $\nu$ on $\Mat_+$ with the 
density 
$$\exp\left\{-\sum_{i,j=1}^n a_{ij} \right\}$$ 
is the probability measure $\mu$ on $\Delta$. In other words, if $A$ is a 
random $n \times n$ matrix with independent standard exponential entries then $X=\Psi(A)$ is 
a random matrix from the simplex $\Delta$ sampled in accordance with the uniform probability measure 
$\mu$. Furthermore, the doubly stochastic scalings of $A$ and $\Psi(A)$ coincide.

Let us choose $A \in \Mat_+$, let $X=\Psi(A)$, and let $B$ be the doubly stochastic scaling 
of $A$. Then 
$$p(X)={N^N \over N!} \per \left(B \otimes J_t \right).$$
In view of Corollary 4.3, the proof of Theorem 3.3 follows from the following result.
\proclaim{(8.1) Proposition} For any $\alpha>0$ there exists $\beta=\beta(\alpha)>0$ such 
that for all positive integers $n$ and $t$ such that 
$$t < e^n$$
the following holds. 

Let $A=\left(a_{ij} \right)$ be the $n \times n$
random matrix with the independent standard exponential entries
and let $B=B(A)$, $B=\left(b_{ij}\right)$, be its doubly stochastic scaling. 
Then 
$$\PP \left\{ \sum_{i=1}^n \left(\max_{j=1, \ldots, n} b_{ij}\right) > \beta \ln N \right\} < N^{-\alpha n}.$$
\endproclaim
\demo{Proof of Proposition 8.1} Let us introduce random variables 
$$\split &u_i={1 \over n} \sum_{j=1}^n \ln  a_{ij} \quad \text{for} \quad i=1, \ldots, n \quad \text{and} \\
&v_j={1\over n}\sum_{i=1}^n \ln  a_{ij} \quad \text{for} \quad j=1, \ldots, n. \endsplit$$

Applying Part 1 of Lemma 7.1, we conclude that for some absolute constant $r >0$ 
we have 
$$\PP \bigl\{u_i \leq -r  \bigr\} \leq 2^{-n} \quad \text{and} \quad \PP\bigl\{v_j \leq -r \bigr\} \leq 2^{-n} \quad 
\text{for all} \quad i,j=1, \ldots, n.$$
It follows then that one can choose a $\beta_1=\beta_1(\alpha)>0$ such that 
$$\aligned &\PP\Bigl\{|i: \ u_i< -r | > \beta_1 \ln N \Bigr\} < {1 \over 6} N^{-\alpha n} \quad \text{and} \\
&\PP \Bigl\{|j: \ v_j < -r |  > \beta_1 \ln N \Bigr\} < {1 \over 6} N^{-\alpha n}. \endaligned \tag8.1.1$$

Applying Part 2 of Lemma 7.1 with $m=n^2$ and using that $t< e^n$ we conclude that for some constant $\beta_2=\beta_2(\alpha)>0$ 
we have 
$$\PP\left\{ {1 \over n^2} \sum_{i,j=1}^n a_{ij} > \beta_2  \right\} < {1 \over 3} N^{-\alpha n}. \tag8.1.2$$

Let us define
$$c_i=\max_{j=1, \ldots, n} a_{ij}.$$
By Part 3 of Lemma 7.1, for some absolute constant $\beta_3=\beta_3(\alpha)>0$ we have 
$$\PP \left\{ {1 \over n} \sum_{i=1}^n c_i > \beta_3 \ln N \right\} < {1 \over 3} N^{-\alpha n}. \tag8.1.3$$ 
Let us define a set ${\Cal A} \subset \Mat_+$ of matrices by 
$$\aligned {\Cal A}=\Biggl\{A=\left(a_{ij}\right): \quad & |i: \ u_i < - r|  \leq \beta_1 \ln N, \quad |j: \ v_j< -r| \leq \beta_1 \ln N, \\ 
& {1 \over n^2} \sum_{i,j=1}^n a_{ij} \leq  \beta_2, \quad \text{and} \\
&{1 \over n} \sum_{i=1}^n c_i \leq \beta_3 \ln N \Biggr\}. \endaligned$$
From (8.1.1)--(8.1.3) we conclude 
$$\PP\bigl\{A \in {\Cal A} \bigr\} \geq 1 - N^{-\alpha n}.$$

Let us pick a matrix $A \in {\Cal A}$ and let $B$ be its doubly stochastic scaling.
By Theorem 6.1, we have 
$$\aligned \ln  b_{ij} \leq {n \over n-2} &\ln a_{ij} - {1 \over n-2} u_i -{1 \over n-2}v_j \\ &+ {n \over n-2} \ln \left({1 \over n} \sum_{i,j=1}^n a_{ij} \right) - {2n-2 \over n-2} \ln(n-1). \endaligned \tag8.1.4$$
Let us define 
$$I=\Bigl\{i: \quad u_i < -r \Bigr\} \quad \text{and} \quad J=\Bigl\{j: \quad v_j < - r \Bigr\},$$
so $|I|, |J| \leq \beta_1 \ln N$.
Thus  from (8.1.2) and (8.1.4) we deduce that for 
some constant $\gamma=\gamma(\beta_1, \beta_2, \beta_3)$ we have
$$b_{ij} \leq {\gamma \over n} a_{ij}^{n/(n-2)}  \leq {\gamma \over n} c_i^{n/(n-2)}  \quad \text{for} \quad i \notin I, \ j \notin J.$$

To complete the proof, we use the estimates
$$\split &\max_{j=1, \ldots, n} b_{ij} \leq 1 \quad \text{for} \quad i \in I \quad \text{and} \\
& \max_{j=1, \ldots, n} b_{ij} \leq {\gamma \over n} c_i^{n/(n-2)} + \sum_{j \in J} b_{ij} \quad \text{for} \quad 
i \notin I. \endsplit$$
Summarizing,
$$\split \sum_{i=1}^n \left( \max_{j=1, \ldots, n} b_{ij} \right) \leq & |I| +  {\gamma \over n} 
\sum_{i=1}^n c_i^{n/(n-2)} 
+ \sum_{i \notin I} \sum_{j \in J} b_{ij} \\ \leq 
&|I| + {\gamma \over n} \left( \sum_{i=1}^n c_i \right)^{n/(n-2)} + |J| \\ \leq &\beta \ln N \endsplit$$
for some $\beta=\beta(\beta_1, \beta_3, \gamma)$
as desired.
{\hfill \hfill \hfill} \qed
\enddemo

\head 9. Proof of Theorem 3.4 \endhead

We use that both $f$ and $\phi$ are {\it positive homogeneous} of degree $N$, that is,
$$\split & f(\lambda X)=\lambda^N f(x) \quad \text{and} \quad \phi(\lambda X)=\lambda^N \phi(X) 
\\ &\quad \text{for all} \quad  X \in \Mat_+(n, n) \quad \text{and} \quad \lambda >0 \endsplit$$ 
and {\it monotone}, that is 
$$\split &f(X) \leq f(Y) \quad \text{and} \quad \phi(X) \leq \phi(Y) \\
&\quad \text{for all} \quad X, Y \in \Mat_+(n, n), \quad X=\left(x_{ij}\right), \quad Y=\left(y_{ij}\right) \\
&\quad \text{such that} \quad x_{ij} \leq y_{ij} \quad \text{for} \quad i,j=1, \ldots, n. \endsplit$$
Among these properties only the monotonicity of $\phi$ is not immediately obvious. It follows, 
for example, from the following representation of $\sigma(X)$, see \cite{MO86} and Section 2.4.
For a positive matrix $X=\left(x_{ij}\right)$ we have  
$$\split n^n\sigma(X)=&\left(\min \sum_{i,j=1}^n x_{ij} \xi_i \eta_j  \right)^n  \quad \text{over all} \quad 
\xi_i, \eta_j >0\\
&\text{subject to} \quad \prod_{i=1}^n \xi_i =\prod_{j=1}^n \eta_j =1.
\endsplit$$ 

Let $dx$ be the Lebesgue measure on the hyperplanes $\sum_{i,j=1}^n x_{ij}=const$ in the 
space of all $n \times n$ matrices $X=\left(x_{ij}\right)$. Using that $f(X)$ is homogeneous, 
we get 
$$\int_{(1-\delta n^2) \Delta} f(x) \ dx = (1-\delta n^2)^{N+n^2-1} \int_{\Delta} f(x) \ dx.$$
On the other hand, for all $X \in \left(1-\delta n^2\right) \Delta$ the matrix $Y=\left(y_{ij}\right)$ defined 
by $y_{ij}=x_{ij} + \delta$ lies in $\Delta_{\delta}$ and $f(Y) \geq f(X)$, which completes the 
proof of Part 1.    

To prove Part 2, let 
$$\alpha=\max_{i,j} \big| x_{ij}-y_{ij}\big|.$$
Hence 
$$\split &x_{ij} \leq y_{ij} + \alpha  \leq \left(1 + \alpha/\delta \right) y_{ij} \quad \text{and} \\
&y_{ij} \leq x_{ij} + \alpha \leq \left(1+ \alpha/\delta \right) x_{ij} \quad \text{for all} \quad i, j \endsplit $$
Using monotonicity and homogeneity of $\phi$ we conclude that 
$$\phi(X) \leq \left(1+\alpha/\delta\right)^N \phi(Y) \quad \text{and} \quad \phi(Y) \leq \left(1+\alpha/\delta\right)^N \phi(X),$$ 
from which the proof follows.
{\hfill \hfill \hfill} \qed

\head 10. Concluding remarks \endhead

\subhead (10.1) Counting general contingency tables \endsubhead
It is plausible to attempt to devise similar algorithms for counting contingency tables 
with the given row and column sums $R=\left(r_1, \ldots, r_m \right)$ and $C=\left(c_1, \ldots, c_n \right)$, where
$$r_1+\ldots +r_m=c_1+ \ldots +c_n=N,$$
 cf. Section 1. While the general idea can be easily generalized to this case,
 cf. \cite{Ba05}, we were unable so far to prove all the necessary bounds, except in 
 the special case when the row sums are equal
 $$r_1 = \ldots = r_m=t $$
or the column sums are equal
$$c_1= \ldots =c_n=t,$$
but not necessarily both.

Suppose, for example, that the row sums are equal. Modifying the construction slightly, one 
can represent the required number of tables by the integral 
$$\int_Q f \ d \mu,$$
where $Q$ is the set of non-negative $m \times n$ matrices with all the row sums equal to 1
(geometrically, $Q$ is a product of $m$ simplices of dimension $(n-1)$ each) and $\mu$ is 
the  Lebesgue measure on $Q$ normalized by the condition $\mu(Q)=1$. The function $f$ factors 
into the product $f=p \phi$ of a log-concave function $\phi$ and a slowly varying function $p$ and 
all the necessary estimate can be carried through, resulting in a randomized polynomial time 
algorithm approximating the number of tables within a factor of $N^{\log N}$ and a randomized 
quasi-polynomial algorithm of $(1/\epsilon)^{O(1)} N^{\log N}$ complexity to approximate the 
number of tables within any given relative error $\epsilon>0$.

\subhead (10.2) Improving the bound \endsubhead The bottleneck of our algorithm is defined 
by the ratio
$$c(n, t)=|\Sigma(n, t)| /\left( \int_{\Delta} \phi \ d \mu \right),$$
where $\phi$ is the log-concave density on the simplex $\Delta$ defined by (2.4.1). Roughly
speaking, $c(n,t)$ is the main contribution to the computational complexity. We proved that
$c(n,t)=N^{O(\ln N)}$ and some conjectural inequalities for the permanent (Conjectures 1.1 and 1.6
of \cite{Sa06}) imply that we can choose the threshold $T=N^{O(1)}$ in Section 1.3 and therefore one should have $c(n,t)=N^{O(1)}$.
 
 \head Acknowledgments \endhead
 
 The authors are grateful to Jes\'us De Loera who computed some of the values 
 of $|\Sigma(n,t)|$ for us using his {\tt LattE} code. The third author would like to thank 
 Radford Neal and Ofer Zeitouni for helpful
discussions.

\enddocument